\documentclass[12pt]{amsart}
\usepackage{amssymb,amsthm,amsmath,amsfonts}
\usepackage{hyperref}
\usepackage{float}
\usepackage[margin=2cm,a4paper]{geometry}
\usepackage[usenames, dvipsnames]{color}
\usepackage{enumerate}
\usepackage[T1]{fontenc}
\usepackage{graphicx}

\numberwithin{equation}{section}
\newtheorem{theorem}{Theorem}[section]
\newtheorem{lemma}[theorem]{Lemma}
\newtheorem{proposition}[theorem]{Proposition}

\theoremstyle{definition}
\newtheorem{example}[theorem]{Example}
\newtheorem{remark}[theorem]{Remark}
\newtheorem{corollary}[theorem]{Corollary}

\newcommand{\N}{\mathbb N}

\DeclareMathOperator{\Tr}{Tr}
\DeclareMathOperator{\tr}{tr}
\newcommand{\Mat}{\mathrm{M}}
\newcommand{\DD}{\mathrm{D}_n}
\newcommand\A{\mathcal{A}}
\newcommand\CC{\mathcal{C}_p}
\newcommand{\XX}{{X}}  
\newcommand{\YY}{Y}
\renewcommand{\SS}{S}
\newcommand{\LL}{K}
\newcommand\KK{L}
\newcommand{\TT}{T}
\newcommand{\MM}{M}
\newcommand{\FF}{F}
\newcommand{\GG}{G}

\newcommand{\ZZ}{Z} 
\newcommand{\Za}{X}
\newcommand{\Zb}{Y}

\newcommand{\MP}{\mathrm{MP}}
\newcommand{\NC}{\mathcal{NC}}
\newcommand{\Krewl}[1]{K{#1}}
\newcommand{\one}{\mathbf{1}} 

\date{}
\begin{document}

\title{On Shepp's Property in Free Probability}

	\author[Wiktor Ejsmont]{Wiktor Ejsmont}
	\address[Wiktor Ejsmont]{ Department of Telecommunications and Teleinformatics, 
		Wrocław University of Science and Technology \\ Wybrze\.ze Wyspia{\'n}skiego 27, 50-370 Wroc\l aw, Poland }
	\email{wiktor.ejsmont@gmail.com}

\author[Takahiro Hasebe]{Takahiro Hasebe}
\address[Takahiro Hasebe]{Department of Mathematics, Hokkaido University\\
 Science Bldg. No. 3, Kita 10, Nishi 8, Kita-Ku, Sapporo, Hokkaido 060-0810, Japan}
 \email{thasebe@math.sci.hokudai.ac.jp}

  \author[Patrycja H\k{e}{\'c}ka-J\k{e}draszczyk]{Patrycja H\k{e}{\'c}ka-J\k{e}draszczyk}
 \address [Patrycja H\k{e}{\'c}ka-J\k{e}draszczyk]{Department of Telecommunications and Teleinformatics, Wrocław University
of Science and Technology\\
Wybrze\.ze Wyspia\'nskiego 27, 50-370 Wroc\l aw, Poland}
\email{patrycja.hecka@gmail.com}

\begin{abstract}
We establish a free analogue of Shepp's theorem for
freely independent standard semicircular elements. 
We show more generally that a family of self-normalized squared quadratic forms is free, and each 
has a standard Marchenko-Pastur distribution. 
The free analogue of Shepp’s
theorem follows by taking the coefficients of quadratic forms to be the three Pauli
  matrices. 
We also prove some new free analogues of Lukacs-type independence results. The key part of the proof is a construction of an equivalent and simpler operator model on a corner algebra compressed by a free projection. 
\end{abstract}
\maketitle
\begin{verse} \textbf{
Dedicated to our dear friend and mentor,
Marek Bożejko,
on the occasion of his 80th birthday.}
\end{verse}
\vskip10pt


\section{Introduction}\label{section:intro}
Free analogues of notions and results in classical probability theory have attracted considerable interest, see e.g.~
\cite{Arizmendi,BP99,BozejkoBryc:2006,EjsFranzSzpoj,EjsLehner,Sakuma,SU25,Szpojankowski2015,Sz17,SzpojankowskiWesolowski:2014}. We are particularly interested in the independence of functions of random variables. An example  in classical probability is Kac-Bernstein's theorem \cite{Bernstein1941,Kac1939}, 
which says that if $X,Y$ are independent random variables, then  $X+Y$ and $X-Y$ are independent if and only if $X, Y$ have normal distributions with the same variance. The free analogue has been shown in \cite[Theorem 5.3]{Nica96}, which characterizes the semicircle distributions with free independence.  
Another example due to Lukacs \cite{Luk55} says that if $\Za$ and $\Zb$ are independent, nondegenerate and positive random variables, then $X/(X+Y)$ and $X+Y$ are independent if and only if $X,Y$ have gamma distributions with the same scale parameter. This also has a free analogue \cite{ BozejkoBryc:2006,Szpojankowski2015}, in which the gamma distributions are replaced with the Marchenko-Pastur distributions. Further results of this kind can be found e.g.~in \cite{Bose,EjsLehner,Yoshida,Sz17}.

In this paper we establish a free analogue of ``squared'' Shepp's theorem. The  classical Shepp's theorem \cite{Shepp} says that 
if $X$ and $Y$ are independent random variables with distribution $\mathrm{N}(0,1)$, then  
\[
\frac{2\Za\Zb}{\sqrt{\Za^2+\Zb^2}}\quad \text{and}\quad \frac{\Za^2-\Zb^2}{\sqrt{\Za^2+\Zb^2}}
\] are also independent and both have distribution $\mathrm{N}(0,1)$. A natural free analogue would be the self-adjoint elements 
\begin{align*}
\LL&:=(\SS^2+\TT^2)^{-\frac1{4}}(\SS\TT+\TT\SS)
          (\SS^2+\TT^2)^{-\frac1{4}}, \notag\\
\KK&:=(\SS^2+\TT^2)^{-\frac1{4}}(\SS^2-\TT^2)
          (\SS^2+\TT^2)^{-\frac1{4}}, 
\end{align*}
where $\SS$ and $\TT$ are freely independent standard semicircular elements. 
However, simulations suggest that $\LL$ and $\KK$ are not free; according to Voiculescu's asymptotic freeness, the three mixed moments in Table \ref{TabMixed}  should vanish in the large $N$ limit if $\LL$ and $\KK$ were free. Also $\LL$ and $\KK$ do not seem to have the distribution of $S$, see Table \ref{Tab1}.

\begin{table}[ht]
\centering
\renewcommand{\arraystretch}{1.5}
\begin{tabular}{|c||c|c|c|c|}
\hline
Mixed moment & $N=100$ & $N=200$ & $N=400$ & $N=800$ \\
\hline
$\tr_N(\mathring{\LL}_N\mathring{\KK}_N)$
& $-0.0017$ & $-2.3\times10^{-5}$ & $-0.0007$ & $-0.0002$ \\

$\tr_N(\mathring{\LL}_N\mathring{\KK}_N
\mathring{\LL}_N\mathring{\KK}_N)$
& $-0.3896$ & $-0.3753$ & $-0.3794$ & $-0.3780$ \\

$\tr_N(\LL_N^2 \KK_N^2)-\tr_N(\LL_N^2)\tr_N(\KK_N^2)$
& $0.1154$ & $0.1115$ & $0.1127$ & $0.1128$ \\
\hline
\end{tabular}
\vspace{0.3cm}
\caption{Numerically estimated mixed moments. $\LL_N, \KK_N$ are sample matrices defined in the same way as $\LL, \KK$ with $(S,T)$ replaced by a simulated pair of independent GUEs of size $N$. $\tr_N$ is the normalized trace $\frac1{N}\Tr_N$ on the matrix algebra $\Mat_N(\mathbb C)$, and 
$\mathring{\LL}_N:=\LL_N - \tr_N(\LL_N)I_N$ and $\mathring{\KK}_N:=\KK_N - \tr_N(\KK_N)I_N$ are the centered elements.}
\label{TabMixed}
\end{table}

\begin{table}[H]
\centering
\begin{tabular}{|c||c|c|c|}
\hline
$k$
& $m_k(\LL_N)$
& $m_k(\KK_N)$ & $m_k(S)$\\
\hline
$1$ & 0.0002 & 0.0025 & 0    \\
$2$ & 0.8707 & 0.8714 & 1    \\
$3$ & 0.0032 & 0.0044 & 0     \\
$4$ & 1.3615 & 1.3676 & 2   \\
$5$ & 0.0124 & 0.0078 & 0   \\
$6$ & 2.5851 & 2.6106 & 5   \\
$7$ & 0.0426 & 0.0157 & 0   \\
$8$ & 5.4368 & 5.5259 & 14  \\
\hline
\end{tabular}
\caption{
Numerically estimated normalized trace moments $m_k(\LL_N)= \tr_N(\LL_N^k)$ and $m_k(\KK_N)= \tr_N(\KK_N^k)$ for matrices $\LL_N, \KK_N$ defined as in Table \ref{TabMixed} with $N=1500$, together with the theoretical moments $m_k(S)$ of the standard semicircular element $\SS$. 
}
\label{Tab1}
\end{table}

We show that the square version of Shepp's theorem turns out to hold in free probability. Observe first that classical Shepp's theorem implies that 
\[
\Gamma_1:=\frac{4\Za^2\Zb^2}{\Za^2+\Zb^2}\quad\text{ and }\quad\Gamma_2:= \frac{(\Za^2-\Zb^2)^2}{\Za^2+\Zb^2} 
\]
are independent and both have the same  distribution as $\Za^2$, which is the gamma distribution $(2\pi x)^{-1/2}e^{-x/2}\one_{(0,\infty)}(x)\, dx$.

\begin{theorem}\label{thm:main1} If $\SS$ and $\TT$ are freely independent standard semicircular elements in a tracial $W^\ast$-probability space $(\A,\varphi)$, then 
\begin{align*}
\MM_0 &:= \SS^2+\TT^2, \\ 
\MM_1&:=(\SS^2+\TT^2)^{-\frac1{2}}(\SS\TT+\TT \SS)^2
          (\SS^2+\TT^2)^{-\frac1{2}}, \notag\\
\MM_2&:=(\SS^2+\TT^2)^{-\frac1{2}}(\SS^2-\TT^2)^2
          (\SS^2+\TT^2)^{-\frac1{2}}, \notag\\
\MM_3&:=(\SS^2+\TT^2)^{-\frac1{2}}\bigl(i(\TT\SS-\SS\TT)\bigr)^2         (\SS^2+\TT^2)^{-\frac1{2}}\notag
\end{align*}
are freely independent, and  $\MM_1, \MM_2,\MM_3$ have the same distribution as $\SS^2$, which is the standard Marchenko-Pastur distribution $(2\pi)^{-1}x^{-1/2}\sqrt{4-x}\,\one_{(0,4)}(x)\,dx$. 
\end{theorem}
\begin{remark}
      The above result for $\MM_1$ and $\MM_2$ is the free analogue of squared Shepp's theorem. The operator $\MM_0$ has the Marchenko-Pastur distribution $\MP(2,1)$, see Section \ref{sec:basics}. In classical probability, there is no counterpart of $\MM_3$. Moreover, $\Gamma_0:=\Za^2+\Zb^2$ is not independent of $\{\Gamma_1,\Gamma_2\}$ as we have the linear relation $\Gamma_0 = \Gamma_1 +\Gamma_2$.   
\end{remark}
   
In fact, we prove a more general result for quadratic forms as follows.   

\begin{theorem}\label{thm:main}Let $\SS_1,\SS_2,\dots,\SS_n$, $n\geq 2$, be freely independent standard semicircular elements in a tracial $W^\ast$-probability space $(\A,\varphi)$. With a matrix 
    $A=[a_{ij}]_{i,j=1}^n\in \Mat_n(\mathbb C)$
    we associate the quadratic form 
    \begin{equation}\label{eq:Q}
       Q(A):=\sum_{i,j=1}^n a_{ij}\SS_i\SS_j.  
    \end{equation}
Let $U_0,U_1,\ldots,U_m \in \Mat_n(\mathbb C)$, $m\ge1,$ be unitary matrices with $U_0=I_n$, and define 
\[
\MM_k=Q(U_0)^{-\frac1{2}}Q(U_k)Q(U_k)^*Q(U_0)^{-\frac1{2}},
\qquad 0\leq k\leq m.
\]
Then the following conditions are equivalent. 
\begin{enumerate}[\rm(1)]

\item\label{item1}
The matrices $U_0,U_1,\ldots,U_m$ are orthonormal
with respect to the normalized Hilbert-Schmidt
inner product:
\[
\operatorname{tr}_n(U_i^*U_j)=\delta_{ij},
\qquad i,j\in\{0,1,\dots, m\},
\]
where $\operatorname{tr}_n :=\frac1{n}\operatorname{Tr}_n$ is the normalized trace. 

\item\label{item2}
The operators $\MM_0,\MM_1,\ldots,\MM_m$ are  freely
independent and $\MM_1,\ldots,\MM_m$ have the standard Marchenko-Pastur
distribution. 

\item\label{item3} The operators $M_1,\ldots,M_m$ have mean one
and are pairwise uncorrelated, that is,
$$
\varphi(M_i)=1, \quad \mathrm{Cov}(M_i,M_j):=\varphi(M_iM_j)
-\varphi(M_i)\varphi(M_j)=0,
\qquad 1\leq i\neq j\leq m. 
$$
The latter covariance condition is unnecessary in case $m=1$. 

\end{enumerate}
\end{theorem}

\begin{remark}
An interesting feature of Theorem~\ref{thm:main}
is that the freeness  condition in \eqref{item2}
follows from the much weaker requirement
of pairwise uncorrelatedness  together with the normalization of  the first moments. 
This phenomenon is reminiscent of the classical
Gaussian setting, where pairwise uncorrelated
components of a jointly Gaussian random vector
are necessarily independent. 
\end{remark}

\begin{remark} The operator  $\MM_0=\SS_1^2+\SS_2^2+\cdots+\SS_n^2$ has the Marchenko-Pastur distribution $\MP(n,1)$, see Section \ref{sec:basics}. 
    Theorem \ref{thm:main} includes Theorem \ref{thm:main1} as the special case of the identity matrix and the three Pauli matrices
\[
I_2=\begin{pmatrix}1&0\\0&1\end{pmatrix},\qquad
U_1=\begin{pmatrix}0&1\\1&0\end{pmatrix},\qquad
U_2=\begin{pmatrix}1&0\\0&-1\end{pmatrix},
\qquad U_3=\begin{pmatrix}0&-i\\i&0\end{pmatrix},
\]
which are orthonormal with respect to the inner product $\langle A,B\rangle:=\tr_2(A^*B)$. 
\end{remark}

The rest of the paper is organized as follows. Section \ref{sec:basics} collects rather basic notions and facts in free probability such as free cumulants and standard semicircular elements. Section \ref{sec:proof} is devoted to the proof of Theorem \ref{thm:main}, in which an equivalent operator model is constructed on a corner algebra compressed by a freely independent projection. In Section \ref{subsec:free-lukacs-corner}, we further apply the techniques of freely compressed corner algebras to free Lukacs-type elements, which illustrates the power of this method.  In Section \ref{sec:alternative} we give an alternative derivation of the Marchenko-Pastur law of the elements $M_1,M_2$ in Theorem \ref{thm:main1}, independently of Section \ref{sec:proof}.

\section{Basics of free probability}\label{sec:basics}

By a tracial $W^*$-probability space we mean a pair $(\mathcal A,\varphi)$,
where $\mathcal A$ is a von Neumann algebra and $\varphi$ is a faithful
normal tracial state. We denote by  $1_\A$ the unit of $\A$. 
If $X\in\A$ is self-adjoint, then the distribution of $X$ with respect to the state $\varphi$ is the unique probability measure $\mu_X$ on $\mathbb{R}$ such that 
\[
\varphi( f(X))=\int_\mathbb{R}f(\lambda)d\mu_X(\lambda)\]
for any bounded Borel function $f$ on $\mathbb{R}$. We also write $X \sim \mu$ to express that $X$ has distribution $\mu$.

Let $I$ be a set of indices. A family  $(\mathcal{A}_i)_{i\in I}$ of $\ast$-subalgebras
of $\mathcal{A}$, each containing $1_\A$, is called \emph{free} or \emph{freely independent} if   
\[
\varphi(X_{1}X_2 \cdots X_{n} ) = 0
\]
whenever $\varphi(X_{j} ) = 0$ for all
$j = 1,\dots,n$ and $X_j\in\mathcal A_{i_j}$ with $i_j\neq i_{j+1}$ for all $j=1,\dots,n-1$.
 A family $(S_i)_{i\in I}$ of subsets $S_{i}\subseteq\A$ is called free (or freely independent) if the $\ast$-subalgebras $\A_i:=\ast$-alg$\{1_\A,S_i\}$ are free. 

 Let $(\mathcal B,\psi)$ be another tracial $W^*$-probability space. Two families $(X_i)_{i\in I} \subseteq\A$ and $(Y_i)_{i\in I} \subseteq\mathcal B$ are said to have the same \emph{joint $\ast$-distribution} if 
 \[
 \varphi(X_{i_1}^{\epsilon_1}\cdots X_{i_n}^{\epsilon_n}) = \psi(Y_{i_1}^{\epsilon_1}\cdots Y_{i_n}^{\epsilon_n})
 \]
for all $n\ge1$,  $i_1,\dots, i_n \in I$ and $\epsilon_1,\dots, \epsilon_n \in \{1,\ast\}$.

We recall some basic notions and results from the
combinatorics of free probability that will be used
in the proof of the main theorem.
For a comprehensive treatment, we refer to
Nica and Speicher \cite{nica}. 
Let $S$ be a totally ordered finite set.
A partition of $S$ is a set of nonempty, mutually disjoint subsets
(called \emph{blocks}) $B_1,B_2,\dots,B_k\subseteq S$ 
whose union is $S$. Any partition $\pi$ defines an equivalence relation on $S$,
denoted by $\sim_\pi$, whose equivalence classes are the blocks of $\pi$. 
That is, $i\sim_\pi j$ if $i$ and $j$ belong to the same block of $\pi$.
A partition $\pi$ is called \emph{noncrossing} 
if different blocks do not interlace, i.e., there is no quadruple of
elements $i<j<k<l$ such that $i\sim_\pi k$ and $j\sim_\pi l$ 
but $i\not\sim_\pi j$. 
The set of noncrossing partitions of $S$ is denoted by $\NC(S)$,
in the case where $S=[n]:=\{1, \dots , n\}$ we write
$\NC(n):=\NC([n])$. 
$\NC(n)$ is a poset under refinement order, where we say $\pi\leq \rho$ if
every block of $\pi$ is contained in a block of $\rho$. 
It turns out that $\NC(n)$ is in fact a lattice,  see \cite[Lecture~9]{nica}. 

Kreweras \cite{Kreweras} discovered an interesting antiisomorphism
of the lattice $\NC(n)$, now called the \emph{Kreweras complementation map},
for which we use the following convention.
Given a noncrossing partition $\pi$ of $\{1,2,\dots,n\}$,
the \emph{Kreweras complement} $\Krewl(\pi)$ is the maximal
noncrossing partition of the ordered set $\{\bar{1},\bar{2},\dots,\bar{n}\}$
such that $\pi\cup\Krewl(\pi)$ is a noncrossing partition of
the interlaced set
$\{1,\bar{1},2,\bar{2},\dots,n,\bar{n}\}$. 
It satisfies $
|\Krewl(\pi)|+|\pi|=n+1
$;
see \cite[Definition 9.21 and Exercise 9.23]
{nica} for more details. 
Given a tracial $W^*$-probability space $(\A,\varphi)$, free cumulants are multilinear maps $\kappa_n \colon \mathcal{A}^n\to\mathbb C$, $n\ge1$,  defined implicitly by the moment-cumulant formula
\begin{align}
\varphi(X_{1}X_{2}\dots X_{n}) &= \sum_{\pi \in \NC(n)}\kappa_{\pi}(X_{1},X_{2},\dots ,X_{n}),\label{eq:DefinicjaKumulant}
\intertext{where}
\kappa_{\pi}(X_{1},X_{2},\dots ,X_{n})&:=\prod_{B \in \pi}\kappa_{|B|}(X_{i}:i \in B). \notag 
\end{align}
We also write $\kappa_n^\A$  to emphasize that $\kappa_n$ depends on $(\A,\varphi)$. 
If the entries are identical, we simply write $\kappa_n(X):=\kappa_n(X,X,\dots, X)$. In particular, $\kappa_1(X)=\varphi(X)$ is the mean and $\kappa_2(X)=\varphi(X^2)-\varphi(X)^2$ is the variance.

Another important tool is a combinatorial
formula for products of freely independent
families. 
 Let $\{X_1,X_2,\dots,X_n\}\subseteq \A$ be a subset free from a subset $\{Y_1,Y_2,\dots,Y_n\}\subseteq \A$. Then
\begin{align}
 \label{productkreweras}
\varphi(X_1Y_1\cdots X_nY_n)
&=
\sum_{\pi\in \NC(n)}
\kappa_\pi(X_1,\ldots,X_n)
\varphi_{K(\pi)}(Y_1,\ldots,Y_n),
\end{align}
where 
\[
\varphi_{K(\pi)}(Y_1,\ldots,Y_n):=\prod_{B \in \Krewl(\pi)} \varphi\left(\overrightarrow{\prod_{i \in B}} Y_i\right) 
\]
and $\overrightarrow{\prod}_{i\in B}Y_i$ is the ordered product, i.e., it is  $Y_{j_1}Y_{j_2}\cdots Y_{j_m}$ when $B=\{j_1,j_2,\dots,j_m\}$ with $j_1 < j_2 < \cdots < j_m$.
This identity follows from
\cite[Chapter~2]{MingoSpeicher}. 

A \emph{standard semicircular element} is a self-adjoint element in $\mathcal A$ that has the \emph{semicircle distribution} 
\[
\frac{\sqrt{4-x^2}}{2\pi}\one_{[-2,2]}(x)\,dx.
\]
A standard semicircular element $X$ has mean zero and  variance one. The higher-order free cumulants $\kappa_n(X)$, $n\ge3$, all vanish.

The \emph{Marchenko-Pastur
(or free Poisson) distribution} $\MP(\lambda,\alpha)$ with parameters 
$\lambda,\alpha>0$ is the probability measure
\begin{equation}\label{MPdist}
\begin{aligned}
\MP(\lambda,\alpha)(dx)
&=(1-\lambda)_+\,\delta_0(dx)
+\frac{\sqrt{(b-x)(x-a)}}{2\pi\alpha x}
\one_{[a,b]}(x)\,dx,\\
a&=\alpha(1-\sqrt\lambda)^2,\qquad
b=\alpha(1+\sqrt\lambda)^2.
\end{aligned}
\end{equation}
In particular, $\MP(1,1)$ is called a \emph{standard Marchenko-Pastur distribution}.  
A self-adjoint element $Y$ has distribution $\MP(\lambda,\alpha)$ if and only if  
\begin{align*}  
\kappa_n(Y)=\lambda \alpha^n,\qquad n \geq 1. 
\end{align*}
If $Y_1, Y_2$ are freely independent self-adjoint elements and $Y_i$ has distribution $\MP(\lambda_i,\alpha)$, then $Y_1+Y_2$ has distribution $\MP(\lambda_1+\lambda_2,\alpha)$ because $\kappa_n(Y_1+Y_2) = \kappa_n(Y_1)+\kappa_n(Y_2)=(\lambda_1+\lambda_2)\alpha^n$. 

Finally, it is well known and easy to check that if $X$ is a standard semicircular element then $X^2$ has distribution $\MP(1,1)$.

\section{Proof of the main theorem}\label{sec:proof}
\subsection{A corner algebra}\label{sectioncorneralgebra}

We introduce a corner algebra used in the main theorem.
Let $0<\lambda<1$ and $\alpha>0$. 
Let $(\operatorname{vN}(\ZZ),\tau_\ZZ)$ be a von Neumann algebra generated by a single element $\ZZ\ge0$ with a normal faithful tracial state $\tau_\ZZ$. We assume that $\ZZ$ has distribution
$\MP(\lambda,\alpha)$.  
Let $n\geq2$ and consider the free product von Neumann algebra
\[
(\mathcal C,\tau)
=
\bigl(\Mat_n(\mathbb C),\tr_n\bigr)
*
\bigl(\operatorname{vN}(\ZZ),\tau_\ZZ\bigr). 
\]
The canonical state $\tau$ is faithful, normal, and tracial,
and the two factors $\Mat_n(\mathbb C)$ and $\operatorname{vN}(\ZZ)$ are freely independent in $(\mathcal C,\tau)$; see \cite[Section~5.3.2]{Popa}. 
Set $p=\mathbf{1}_{(0,\infty)}(\ZZ)$. Since
$\mu_\ZZ(\{0\})=1-\lambda$, we have $\tau(p)=\lambda$.
We write
\[
\CC:=p\mathcal Cp,
\qquad
\tau_p(a):=\frac{\tau(a)}{\lambda},
\quad a\in\CC.
\]
Thus $\CC$ is a von Neumann algebra
with unit $p$ and faithful normal tracial state $\tau_p$.
By the faithfulness of $\tau$, the spectrum of $\ZZ$ in $\mathcal C$
coincides with the support of its distribution:
\[
\sigma_{\mathcal C}(\ZZ)=\{0\}\cup[a_-,a_+],
\qquad
a_\pm=\alpha(1\pm\sqrt{\lambda})^2.
\]
As $p\ZZ^\beta=\ZZ^\beta p =\ZZ^\beta$ for $\beta>0$, we have $\ZZ^\beta \in \CC$. Since $\sigma_{\CC}(\ZZ)=[a_-,a_+]$ and $a_->0$, $\ZZ^\beta$ is invertible in $\CC$.
We denote by $\ZZ^{-\beta} \in \CC$ the inverse of $\ZZ^{\beta}$ for $\beta>0$.

\subsection{Complementary facts}
Below, we recall several results from the literature and prove auxiliary lemmas needed for the main theorem.
In this section $(\A,\varphi)$ always denotes a tracial $W^\ast$-probability space.

\begin{lemma}[{\cite[Proposition 4.5, Theorem 4.7] {EjsLehner}} or {\cite[Proposition 2.2 in a special case]{Yoshida}}] \label{prop:CykliczneVariancja}
  Let $Q(A)\in(\A,\varphi)$ be the quadratic form defined by \eqref{eq:Q}. Let $A_1,A_2,\dots, A_m \in \Mat_n(\mathbb C)$, $m \ge1$. 
  Then 
    \begin{equation}  \label{eq:kumulantsamplevariance}
      \kappa_m^\A(Q(A_1),Q(A_2),\dots,Q(A_m))= 
      \Tr_n( A_1A_2\cdots A_m).
    \end{equation} 
 \end{lemma}

\begin{proposition}\label{lem:rownoscrozkladow}
Let $Q(A) \in (\A,\varphi)$ be the quadratic form defined by \eqref{eq:Q}. 
In the setting of Section \ref{sectioncorneralgebra}, assume that $\ZZ$ has distribution $\MP(
 1/n,n)$ with $n\ge2$. 
For $A \in \Mat_n(\mathbb C)$ we define
\begin{equation}\label{eq:tQ}
    \widetilde Q(A)=\ZZ^{\frac1{2}}A\ZZ^{\frac1{2}}\in \CC,\quad
 A \in \Mat_n(\mathbb C).
\end{equation}

Let $A_1,A_2,\ldots,A_m \in \Mat_n(\mathbb C)$, $m\ge1$. Then the families  \[
\bigl( Q(A_1),\ldots,Q(A_m)\bigr)\in (\A,\varphi)
\quad\text{and}\quad
\bigl( \widetilde Q(A_1),\ldots,\widetilde Q(A_m)\bigr)\in (\CC,\tau_p)
\] have the same joint
$\ast$-distribution.

\end{proposition}
\begin{proof} 
Traciality and the product formula~\eqref{productkreweras} yield
\[
\tau_p\bigl(\widetilde Q(A_1)\cdots\widetilde Q(A_m)\bigr)
=n\tau(\ZZ A_1\ZZ A_2\cdots \ZZ A_m)
=n\sum_{\pi\in\NC(m)}\kappa_\pi^{\mathcal C}(\ZZ)\,
(\tr_{n})_{K(\pi)}(A_1,\ldots,A_m).
\]
Here $\kappa_m^{\mathcal C}(\ZZ)=\frac1n\,n^m=n^{m-1}$ and so  $\kappa_\pi^{\mathcal C}(\ZZ)=n^{m-|\pi|}$, while
$(\tr_n)_{K(\pi)}$ is the product of normalized traces of the cyclic
matrix products over the blocks of $K(\pi)$. 
Since $|K(\pi)|=m+1-|\pi|$, all powers of $n$ cancel, leaving
\[
\tau_p\bigl(\widetilde Q(A_1)\cdots\widetilde Q(A_m)\bigr)
=\sum_{\pi\in\NC(m)}
\prod_{B\in K(\pi)}
\Tr_n\left(\prod_{j\in B}^{\longrightarrow}A_j\right).  
\]
Because the Kreweras complementation map $\Krewl{}$ is a bijection of $\NC(m)$ onto itself, using  Lemma~\ref{prop:CykliczneVariancja} we get 
\begin{align*}
\tau_p\bigl(\widetilde Q(A_1)\cdots\widetilde Q(A_m)\bigr) 
&=\sum_{\rho\in\NC(m)}
\prod_{B\in \rho}
\Tr_n\left(\prod_{j\in B}^{\longrightarrow}A_j\right)\\
&=\sum_{\rho\in\NC(m)}
\prod_{B\in\rho}
\kappa_{|B|}^\A(Q(A_i): i\in B)
=\varphi(Q(A_1)\cdots Q(A_m)).
\end{align*} 
The equality also covers moments involving adjoints, since
$Q(A)^*=Q(A^*)$ and
$\widetilde Q(A)^*=\widetilde Q(A^*)$.
\end{proof}

\begin{lemma}\label{lem:conjugates}
Let $\mathcal B$ be a $\ast$-subalgebra of  $(\A,\varphi)$ containing  $1_\A$, and let $u_1,\ldots,u_m\in\A$, $m\ge1$, be unitary
elements such that $\mathcal B$ is free from the set $\{u_1,u_2,\dots,u_m\}$, and  $\varphi(u_i^*u_j)=\delta_{ij}$ for all $i,j\in\{1,2,\dots,m\}$.
Then the subalgebras $u_i\mathcal B u_i^*$, $1\le i \le m$, are free.
\end{lemma}
\begin{proof}
Take centered elements $u_ib_iu_i^\ast\in u_i\mathcal B u_i^*$ for $i=1,2,\dots,m$, i.e.,   $\varphi(b_i)=0$, and indices $i_1,\dots,i_n$ with $i_k\neq i_{k+1}$ for $k=1,\dots,n-1$. By traciality,
	\[
		\begin{aligned}
			&\varphi\bigl((u_{i_1}b_{i_1}u_{i_1}^\ast)(u_{i_2}b_{i_2}u_{i_2}^\ast)\cdots
			(u_{i_n}b_{i_n}u_{i_n}^*)\bigr)=\varphi\bigl(b_{i_1}(u_{i_1}^\ast u_{i_2})b_{i_2}(u_{i_2}^\ast u_{i_3})\cdots
			b_{i_n}(u_{i_n}^*u_{i_1})\bigr).
		\end{aligned}
		\]
Each factor $u_{i_k}^*u_{i_{k+1}}$ with $k<n$ is centered.
The final factor is also centered unless $i_n=i_1$; in that case
it is the identity and may be omitted.
We are therefore left with an alternating product of free centered
elements,
whose trace is zero.
\end{proof}

\begin{lemma}[Free compression]\label{lem:compression}

Consider elements $Y_1,\ldots,Y_m\in\mathcal{A}$ and a $\ast$-subalgebra $\mathcal B$ containing $1_\A$ such that $\mathcal B$ is free from $\{Y_1,\ldots,Y_m\}$. Let 
$p\in \mathcal B$ be a projection with $\varphi(p)\neq 0$.  
Let $\varphi_p(a) := \frac{\varphi(a)}{\varphi(p)}$ for $a \in p \A p$. 
Let $\kappa_n^{p\A p}$ be the free cumulants in the compressed tracial $W^\ast$-probability space $(p\A p, \varphi_p)$. Then   
\begin{equation}
\kappa_n^{p\A p}
\bigl(pY_{i_1}p,\ldots,pY_{i_n}p\bigr)
=
\frac{1}{\lambda}
\kappa_n^\A
\bigl(\lambda Y_{i_1},\ldots,\lambda Y_{i_n}\bigr),
\label{eq:compression-formula}
\end{equation}
for all $n\geq 1$ and all
$i_1,\ldots,i_n\in\{1,2,\dots, m\}$, where 
$
\lambda:=\varphi(p).$ Moreover, suppose that $Y_1, Y_2,\dots, Y_m$ are free in $(\A,\varphi)$. Then for any $X \in \mathcal B$, the elements $pXp, pY_1p, pY_2 p , \dots, p Y_m p$ are free in $(p\A p, \varphi_p)$. 
\end{lemma}
\begin{proof} Formula \eqref{eq:compression-formula} is due to \cite[Application~1.11]{NS96} or \cite[Theorem~14.10]{nica}. According to \eqref{eq:compression-formula}, if $Y_1,\dots, Y_m$ are free, then the mixed cumulants of $pY_1p, pY_2 p , \dots, p Y_m p$ vanish in the compressed space, so that they are free in $(p\A p, \varphi_p)$ (\cite[Corollary 1.12]{NS96} or \cite[Corollary 14.12]{nica}). Moreover,  \cite[Application~1.13]{NS96} implies that $p\mathcal B p$ is free from $\{pY_1p, pY_2 p , \dots, p Y_m p\}$ in $(p\A p, \varphi_p)$. Combining these two freeness facts, together with the associativity, yields the  freeness of $pXp, pY_1p, pY_2 p , \dots, p Y_m p$. 
\end{proof}

\subsection{Proof of Theorem \ref{thm:main}} 
We are now ready to prove our main result.
Since $\SS_1^2,\SS_2^2, \ldots,\SS_n^2$ are free and each has
distribution $\MP(1,1)$, their sum $Q(U_0)$ has distribution
$\MP(n,1)$. By the faithfulness of $\varphi$,
\[
\sigma(Q(U_0))=[a,b],
\qquad
a=(\sqrt n-1)^2,\quad b=(\sqrt n+1)^2.
\]
As $n\geq2$, we have $a>0$, so $Q(U_0)^{-1/2}$ is bounded.
Let $(\mathcal C,\tau)$  be the von Neumann algebra  and $(\CC,\tau_p)$ be the corresponding  corner
constructed in Section \ref{sectioncorneralgebra}, where $\ZZ\sim\MP(\frac1n,n)$. Recall that 
\[
p=\one_{(0,\infty)}(\ZZ),\qquad
\tau(p)=\frac1n,\qquad
\tau_p=n\tau\big|_{\CC}. 
\]
Under the notation of \eqref{eq:tQ},  we have  $\ZZ=\widetilde{Q}(U_0)$.  
By Proposition~\ref{lem:rownoscrozkladow}, the families
\[
\bigl(Q(U_0),Q(U_1),\ldots,Q(U_m)\bigr)
\quad\text{and}\quad
\bigl(\ZZ,\widetilde Q(U_1),\ldots,\widetilde Q(U_m)\bigr)
\]
have the same joint $*$-distribution with respect to
$\varphi$ and $\tau_p$, respectively. 
The unit of $\CC$ is $p$, and the
spectrum of $\ZZ$ in this corner is contained in $[a,b]$.
The function $f(t)=t^{-1/2}$ is continuous on $[a,b]$.
By the Weierstrass approximation theorem, there exist
polynomials $f_N$~($N\in\N$) such that $f_N(t) \to t^{-1/2}$ uniformly on $[a,b]$. 
Continuous functional calculus therefore yields
\[
f_N(Q(U_0))\longrightarrow Q(U_0)^{-\frac1{2}},
\qquad
f_N(\ZZ)\longrightarrow \ZZ^{-\frac1{2}}
\]
in operator norm, with the second convergence
understood in the corner algebra. Norm-continuity
of the states permits passage to the limit in all
joint $*$-moments containing these operators. 
Consequently, $(\MM_0,\MM_1,\ldots,\MM_m)$ in $(\A,\varphi)$ has the same joint
$*$-distribution as $(\widetilde \MM_0,\widetilde \MM_1,\ldots,\widetilde \MM_m)$ in 
$(\CC,\tau_p)$,
where 
\[
\begin{aligned}
\widetilde \MM_k
&:=\ZZ^{-\frac1{2}}\widetilde Q(U_k)
  \widetilde Q(U_k)^*\ZZ^{-\frac1{2}}\\
&=\ZZ^{-\frac1{2}}
  (\ZZ^{\frac1{2}}U_k\ZZ^{\frac1{2}})
  (\ZZ^{\frac1{2}}U_k^*\ZZ^{\frac1{2}})
  \ZZ^{-\frac1{2}}\\
&=pU_k\ZZ U_k^*p, \qquad k\in\{0,1,\ldots,m\}.
\end{aligned}
\]
Thus
\begin{equation}\label{eq:stardist}
(\MM_0,\MM_1,\ldots,\MM_m)
\overset{*\text{-dist}}{=}
\bigl(\ZZ,pU_1\ZZ U_1^*p,\ldots,pU_m \ZZ U_m^*p\bigr).
\end{equation}

\medskip
\noindent
\textbf{Proof of $ \eqref{item1}\Rightarrow\eqref{item2}$}.
We apply Lemma~\ref{lem:conjugates} to the unitaries
$U_0=I_n,U_1,\ldots,U_m$. By the assumed orthonormality,
the algebras
\[
\operatorname{vN}(\ZZ),\quad
U_1\operatorname{vN}(\ZZ)U_1^*,\quad\ldots,\quad
U_m\operatorname{vN}(\ZZ)U_m^*
\]
are free in $(\mathcal C,\tau)$. Set
\[
\ZZ_k=U_k\ZZ U_k^*,\qquad k=1,\ldots,m.
\]
Then $\ZZ_1,\ldots,\ZZ_m$ are free, and
$\operatorname{vN}(\ZZ)$ is free from $\{\ZZ_1,\ldots,\ZZ_m\}$. 
Applying Lemma~\ref{lem:compression} with
$\lambda=\tau(p)=1/n$, 
we have that $\ZZ=p\ZZ p,p\ZZ_1p,\ldots,p\ZZ_mp$ are
 free in the corner algebra $\CC$. By \eqref{eq:stardist},
the same holds for $\MM_0,\MM_1,\ldots,\MM_m$. 
Finally, unitary conjugation preserves the distribution
with respect to $\tau$, so each $\ZZ_k$ has distribution
$\MP(1/n,n)$ in $(\mathcal C,\tau)$. Therefore, for every $r\geq1$, formula \eqref{eq:compression-formula} yields 
\[
\kappa_r^{\CC}(p\ZZ_kp)
=
n^{1-r}\kappa_r^{\mathcal C}(\ZZ_k)
=
n^{1-r}n^{r-1}
=1, \qquad 1 \le k \le m. 
\]
These are precisely the free cumulants of $\MP(1,1)$.
Thus each $\MM_k$ has distribution $\MP(1,1)$ for all $1\le k \le m$.

\medskip
\noindent
 $\eqref{item2}\Rightarrow\eqref{item3}$ is obvious. 

\medskip
\noindent
\textbf{Proof of $\eqref{item3}\Rightarrow\eqref{item1}$}. It suffices to show that $U_i$ has vanishing trace for $1 \leq i \leq m$ and that, if $m\ge2$,
$\operatorname{tr}_n(U_i^*U_j)=0$
for every pair of distinct indices $i,j\in\{1,2,\dots,m\}$.

We first compute the first moment of a 
unitary $U:=U_i$. Write
$
\widetilde \MM_U=pU\ZZ U^*p 
$, which has mean $1$ in $(\CC,\tau_p)$ by assumption \eqref{item3} and the joint distributional
identity \eqref{eq:stardist}. 
Since $p\ZZ=\ZZ$, $\tau(\ZZ)=(1/n)\times n=1$,   $\tau(p)=1/n$,
freeness of $\operatorname{vN}(\ZZ)$ from $\Mat_n(\mathbb C)$ gives
\[
\begin{aligned}
\tau(\widetilde \MM_U)=\tau(pU\ZZ U^*)
&=\tau(p)\tau(\ZZ)\tau(UU^*)+
\bigl(\tau(p\ZZ)-\tau(p)\tau(\ZZ)\bigr)
\tau(U)\tau(U^*)\\
&=\frac1{n}+\left(1-\frac1{n}\right)|\operatorname{tr}_n(U)|^2.
\end{aligned}
\]
Consequently,
$
\tau_p(\widetilde \MM_U)
=
1+(n-1)|\operatorname{tr}_n(U)|^2.
$
Therefore, $\operatorname{tr}_n(U_i)=0$
for every $i=1,\ldots,m$.

We next compute the covariance of two
unitaries $U=U_i$ and $V:=U_j$ with $i\ne j$ in case $m\ge2$.  
Set 
$\mathring{\ZZ}:=\ZZ-1_\mathcal{C},$ and $q:=p-(1/n) 1_\mathcal{C}.$
Then
$
\tau(\mathring{\ZZ})=\tau(q)=0,$ 
$\tau(\mathring{\ZZ}^2)=\tau(\ZZ^2-2\ZZ+1_\mathcal{C})=n+1-2+1=n$,
and
\[
\widetilde \MM_U-p=pU\mathring{\ZZ}U^*p,\qquad
\widetilde \MM_V-p=pV\mathring{\ZZ}V^*p.
\]
By traciality and $p^2=p$,
\[
\mathrm{Cov}(\widetilde \MM_U,\widetilde \MM_V):=\tau_p\bigl((\widetilde \MM_U-p)(\widetilde \MM_V-p)\bigr)
=
n\,\tau(pU\mathring{\ZZ}U^*pV\mathring{\ZZ}V^*).
\]
Expanding both occurrences of $p=(1/n) 1_\mathcal{C}+q$, we obtain
\[
\tau(pU\mathring{\ZZ}U^*pV\mathring{\ZZ}V^*)
=
\frac1{n}|\operatorname{tr}_n(U^*V)|^2.
\]
Indeed, the term containing two occurrences of $q$
is an alternating product of centered elements
from the two free algebras. The terms containing
exactly one occurrence of $q$ also vanish from the traciality and freeness. The term containing no $q$ can be computed as 
\[
\begin{aligned}
\tau(U\mathring{\ZZ}U^*V\mathring{\ZZ}V^*)
&=\tau\bigl(\mathring{\ZZ}(U^*V)\mathring{\ZZ}(V^*U)\bigr)\\
&=\tau(\mathring{\ZZ}^2)\,
  \operatorname{tr}_n(U^*V)
  \operatorname{tr}_n(V^*U)\\
&=n|\operatorname{tr}_n(U^*V)|^2.
\end{aligned}
\]
Thus, 
\[
\tau_p\bigl((\widetilde \MM_U-p)(\widetilde \MM_V-p)\bigr)
=
|\operatorname{tr}_n(U^*V)|^2.
\]
The covariance assumption in \eqref{item3} together with \eqref{eq:stardist} implies that the covariance of $\widetilde \MM_U$ and $\widetilde \MM_V$ also vanishes. We can therefore obtain  
\[
0
=
\tau_p\bigl((\widetilde \MM_U-p)(\widetilde \MM_V-p)\bigr)
=
|\operatorname{tr}_n(U^*V)|^2,
\]
and hence
$
\operatorname{tr}_n(U^\ast V)=0.
$
This finishes the proof of Theorem \ref{thm:main}.

\begin{remark} The assumption that the unitary matrices have vanishing first moments, i.e., $
\operatorname{tr}_n(U_k)
=\operatorname{tr}_n(U_0^*U_k)
=\delta_{0k}=0$ for $k=1,2,\dots,m$, is essential. 
Indeed, let $n=2$ and consider the matrices
\[
U_1=I_2=
\begin{pmatrix}
1 & 0\\
0 & 1
\end{pmatrix},
\qquad
U_2=
\begin{pmatrix}
1 & 0\\
0 & -1
\end{pmatrix}.
\]
Both matrices are self-adjoint and unitary. Moreover,
$
\Tr_2(U_1U_2)=\Tr_2(U_2)=0 \text{ and }\Tr_2(U_1)=2.
$
However, in this case,
$
\MM_1=Q(U_0)^{-1/2}Q(U_1)^2Q(U_0)^{-1/2}
=(\SS_1^2+\SS_2^2)^{-1/2}
  (\SS_1^2+\SS_2^2)^2
  (\SS_1^2+\SS_2^2)^{-1/2} 
 =\SS_1^2+\SS_2^2
\sim\MP(2,1),
$
rather than $\MP(1,1)$. Similar counterexamples can be constructed in higher dimensions.   
\end{remark}

\begin{example}[Weyl-Heisenberg matrices] In addition to the example of Pauli matrices given in introduction,  we present an interesting example involving
the clock and shift matrices, which generate the
finite Weyl-Heisenberg group; see for example
\cite{Farenick,Shahbeigi}. Their products form a family of $n^2$ pairwise orthogonal unitary matrices.
Let $n\geq 2$ and let $\omega=e^{2\pi i/n}$.
Consider the clock matrix
\[
C=\operatorname{diag}
(1,\omega,\omega^2,\ldots,\omega^{n-1})
\]
and the shift matrix
\[
S=
\begin{pmatrix}
0&0&\cdots&0&1\\
1&0&\cdots&0&0\\
0&1&\cdots&0&0\\
\vdots&&\ddots&&\vdots\\
0&0&\cdots&1&0
\end{pmatrix}.
\]
These matrices are unitary and satisfy $
C^n=S^n=I_n,
\text{ }
CS=\omega SC.$
They generate the finite Weyl-Heisenberg group
$
H_n=
\{\omega^r C^k S^\ell:
0\leq r,k,\ell<n\}.
$ Now consider the family of $n^2$ unitary matrices
\[
U_{k,\ell}=C^kS^\ell,
\qquad 0\leq k,\ell<n.
\]
These matrices are pairwise orthogonal with 
\[
\operatorname{Tr}_n
\bigl(U_{k,\ell}U_{p,q}^*\bigr)
=n\delta_{kp}\delta_{\ell q}.
\]
Moreover,
\[
\operatorname{Tr}_n(U_{k,\ell})=0,
\qquad (k,\ell)\neq(0,0).
\]
Consequently, the family
\[
\{U_{k,\ell}:0\leq k,\ell<n\}
\]
forms an orthogonal basis of $\Mat_n(\mathbb C)$,
consisting entirely of unitary matrices.
In particular, this construction provides
$n^2-1$ pairwise orthogonal, traceless unitary
matrices in $\Mat_n(\mathbb C)$, from which we can construct a sufficiently large family of random variables satisfying the assumptions of Theorem \ref{thm:main}. 
\end{example}

\section{The free Lukacs property }
\label{subsec:free-lukacs-corner}

In this section, we establish other applications of
the corner algebra construction introduced above. The results concern the independence of Lukacs-type expressions, though we believe that our method may have more applications. Some details of the proof are omitted, as they have already been used in the preceding calculations.
As before, $(\A,\varphi)$ is a tracial $W^\ast$-probability space throughout this section. 

\begin{proposition}\label{thm:main2}Let $\SS_1,\SS_2,\dots,\SS_n$, $n\geq 2$, be freely independent standard semicircular elements in $(\A,\varphi)$, and $Q(A)$ be the quadratic form in \eqref{eq:Q}. 
Let $A_0=I_n$, and let $A_1,\ldots,A_m\in\Mat_n(\mathbb C)$,
$m\geq 1$, be arbitrary matrices. We define
\[
\MM_0:=Q(A_0),\qquad
\MM_k:=Q(A_0)^{-1/2}Q(A_k)Q(A_0)^{-1/2},
\quad 1\leq k\leq m.
\]
Then $\MM_0$ is free from
 $\{\MM_1,\ldots,\MM_m\}$.
Moreover, for every $r\geq 1$, every
$i_1,\ldots,i_r\in\{1,\ldots,m\}$, and every
$\varepsilon_1,\ldots,\varepsilon_r\in\{1,*\}$,
\begin{equation}\label{eq:angular-cumulants}
\kappa_r^{\A}\bigl(
\MM_{i_1}^{\varepsilon_1},\ldots,
\MM_{i_r}^{\varepsilon_r}\bigr)
=
\frac{1}{n^{r-1}}\,
\kappa_r^{\mathrm{M}_n(\mathbb C)}\bigl(
A_{i_1}^{\varepsilon_1},\ldots,
A_{i_r}^{\varepsilon_r}\bigr).
\end{equation}
\end{proposition}
\begin{proof}
Using exactly the same notation and techniques as at the beginning of the proof of Theorem~\ref{thm:main}, we obtain
\begin{equation}\label{eq:stardist2}
(M_0,\MM_1,\ldots,\MM_m)\text{ in }(\A,\varphi)
\overset{*\text{-dist}}{=}
(Z,pA_1p,\ldots,pA_mp)
\text{ in }(\CC,\tau_p).
\end{equation}
Indeed,
\[
\ZZ^{-1/2}\widetilde Q(A)\ZZ^{-1/2}
=\ZZ^{-1/2}\ZZ^{1/2}A\ZZ^{1/2}\ZZ^{-1/2}
=pAp.
\]
Now $p\in\operatorname{vN}(\ZZ)$, and the unital $*$-algebra
$\Mat_n(\mathbb C)$ is free from $\operatorname{vN}(\ZZ)$
in $(\mathcal C,\tau)$. By
\cite[Application~1.13]{NS96}, the $\ast$-algebra generated by
$pA_1p,\ldots,pA_mp$ is free from
$p\operatorname{vN}(\ZZ)p$ in $(\CC,\tau_p)$.
Since $\ZZ=p\ZZ p\in p\operatorname{vN}(\ZZ)p$,
Equation~\eqref{eq:stardist2} proves that $\MM_0$ is free
from the  $\ast$-algebra generated by $\MM_1,\ldots,\MM_m$.
Finally, the free-compression cumulant formula with
$\lambda=1/n$ gives
\[
\kappa_r^{\CC}(pA_{i_1}p,\ldots,pA_{i_r}p)
=\frac{1}{n^{r-1}}
\kappa_r^{\mathrm{M}_n(\mathbb C)}(A_{i_1},\ldots,A_{i_r}),
\]
and we obtain~\eqref{eq:angular-cumulants}.
\end{proof}

\begin{corollary}\label{cor:wniosek1}
Let $X$ and $Y$ be freely independent self-adjoint elements in $(\A,\varphi)$
with common distribution $\MP(1,1)$. Then the operators
\[
X+Y
\quad\text{and}\quad
(X+Y)^{-\frac1{2}}X(X+Y)^{-\frac1{2}}
\]
are freely independent.
\end{corollary}
\begin{proof}
Let $\SS_1$ and $\SS_2$ be freely independent standard semicircular elements. Since $\SS_1^2$ and $\SS_2^2$ are free and both have
distribution $\MP(1,1)$, the pairs
$(X,Y)$ and $(\SS_1^2,\SS_2^2)$ have the
same joint distribution. 
Applying Proposition ~\ref{thm:main2} with $n=2$ and
\[
A=
\begin{pmatrix}
1&0\\
0&0
\end{pmatrix},
\]
we conclude that
\[
\SS_1^2+\SS_2^2
\quad\text{and}\quad
(\SS_1^2+\SS_2^2)^{-\frac1{2}}
\SS_1^2
(\SS_1^2+\SS_2^2)^{-\frac1{2}}
\]
are free.
\end{proof}

Corollary \ref{cor:wniosek1} is a special case of the result established
in \cite{Szpojankowski2015}, i.e., Lemma \ref{szpojan} below. However, the techniques
developed in the present paper allow us to extend
Lemma \ref{szpojan} to a multivariate setting, as shown below.
This the free analogue of the classical counterpart  that the sum of independent gamma
random variables with the same scale parameter is independent of their
normalized proportions, which can be easily verified by computing the joint density function. 

\begin{theorem}
    [Multivariate free Lukacs property]
	\label{prop:multivariate-lukacs}
	Let $X_1,\ldots,X_n$, $n\geq2$, be freely
	independent positive operators in $(\A,\varphi)$ such that
	\[
	X_j\sim\MP(\lambda_j,\alpha),
	\qquad 1\leq j\leq n,
	\]
	where $\lambda_j>0$, $\alpha>0$, and
	$
	\lambda:=\sum_{j=1}^n\lambda_j>1.
	$
	Set
	\[
	F=\sum_{j=1}^nX_j,
	\qquad
	B_j=F^{-\frac1{2}}X_jF^{-\frac1{2}}.
	\]
	Then $F$ is free from 
	 $\{B_1,\ldots,B_n\}$.
 \end{theorem}
\begin{proof}
	Let $(\DD,\rho)$ be the unital $\ast$-algebra of 
	diagonal matrices of size $n$ with the weighted trace
	\[
	\rho(D):=\sum_{j=1}^n\frac{\lambda_j}{\lambda}d_j,\qquad D=\mathrm{diag}(d_1,d_2,\dots,d_n). 
	\]
Then $\rho$ is a normal faithful tracial state on
$\DD$. 
	Consider the free product
	\[
	(\mathcal C,\tau)
	=
	(\DD,\rho)
	*
	(\operatorname{vN}(Z),\tau_Z),
	\]
	where $Z$ has distribution $\MP(1/\lambda,\lambda\alpha)$ in $(\operatorname{vN}(Z),\tau_Z)$. 
	We consider the corner algebra
	$\CC$, with
	$
	p=\mathbf{1}_{(0,\infty)}(Z),
	$
	equipped with the normalized trace
	\[
	\tau_p(X)=\lambda\tau(X),
	\qquad X\in\CC.
	\]
	For $j=1,\ldots,n$, let $E_j\in \DD$ be the orthogonal projection whose $j$th diagonal entry is 1 and the other entries are 0. We  define
	\[
	\widetilde X_j=Z^{\frac1{2}}E_jZ^{\frac1{2}}.
	\]
	Fix $r\geq1$ and indices
	$i_1,\ldots,i_r\in\{1,\ldots,n\}$.
	By traciality and the definition of
	$\tau_p$, we have
	\begin{align}\label{eq:lukacs-kreweras}
		\tau_p\bigl(
		\widetilde X_{i_1}\cdots\widetilde X_{i_r}
		\bigr)
		&=\lambda\tau\bigl(
		Z E_{i_1}Z E_{i_2}\cdots Z E_{i_r}
		\bigr) \notag\\
		&=\lambda\sum_{\pi\in\NC(r)}
		\kappa_\pi^\mathcal{C}(Z)\,
		\rho_{K(\pi)}
		(E_{i_1},\ldots,E_{i_r}).
	\end{align}
	Since $Z\sim\MP(1/\lambda,\lambda\alpha)$, we obtain
	$
	\kappa_s^\mathcal{C}(Z)
	=
	\frac1\lambda(\lambda\alpha)^s
	=
	\lambda^{s-1}\alpha^s.
	$
	Therefore, for every $\pi\in\NC(r)$,
	\[
	\begin{aligned}
		\kappa_\pi^\mathcal{C}(Z)
		&=
		\prod_{V\in\pi}
		\lambda^{|V|-1}\alpha^{|V|}=
		\lambda^{r-|\pi|}\alpha^r.
	\end{aligned}
	\]
	Substituting this into
	\eqref{eq:lukacs-kreweras}, we get
	\begin{align*}
		\tau_p\bigl(
		\widetilde X_{i_1}\cdots
		\widetilde X_{i_r}\bigr)
		&=
		\lambda\sum_{\pi\in\NC(r)}
		\lambda^{r-|\pi|}\alpha^r\,
		\rho_{K(\pi)}
		(E_{i_1},\ldots,E_{i_r})\\
		&=
		\alpha^r\sum_{\pi\in\NC(r)}
		\lambda^{r+1-|\pi|}
		\rho_{K(\pi)}
		(E_{i_1},\ldots,E_{i_r}).
		\intertext{By the fundamental identity
			$|\pi|+|K(\pi)|=r+1$, we obtain}
		&=
		\alpha^r
		\sum_{\pi\in\NC(r)}
		\lambda^{|K(\pi)|}
		\rho_{K(\pi)}
		(E_{i_1},\ldots,E_{i_r}).
		\intertext{The Kreweras complementation map
			is a bijection of $\NC(r)$ onto itself.
			Making the change of variables
			$\sigma=K(\pi)$, we obtain}
		&=
		\alpha^r\sum_{\sigma\in\NC(r)}
		\lambda^{|\sigma|}
		\prod_{V\in\sigma}
		\rho\left(
		\prod_{j\in V}^{\longrightarrow}E_{i_j}
		\right)\\
		&=
		\sum_{\sigma\in\NC(r)}
		\prod_{V\in\sigma}
		\left[
		\lambda\alpha^{|V|}
		\rho\left(
		\prod_{j\in V}^{\longrightarrow}E_{i_j}
		\right)
		\right].
	\end{align*}
	Using the uniqueness of free
	cumulants, we conclude that
	\begin{equation}\label{eq:lukacs-cumulants}
		\begin{aligned}
			&\kappa_r^{\CC}
			(\widetilde X_{i_1},\widetilde X_{i_2},\ldots,
			\widetilde X_{i_r})=
			\lambda\alpha^r
			\rho(E_{i_1}E_{i_2}\cdots E_{i_r}).
		\end{aligned}
	\end{equation}
	Since the projections $E_1,\ldots,E_n$
	are mutually orthogonal, we have
	\[
	E_{i_1}E_{i_2}\cdots E_{i_r}
	=
	\begin{cases}
		E_j, & i_1=\cdots=i_r=j,\\
		0, & \text{otherwise}.
	\end{cases}
	\]
	Consequently, since $\rho(E_i)=\lambda_i/\lambda$ we obtain 
	\[
	\begin{aligned}
		&\kappa_r^{\CC}
		(\widetilde X_{i_1},\ldots,
		\widetilde X_{i_r})=
		\begin{cases}
			\lambda_j\alpha^r,
			& i_1=\cdots=i_r=j,\\
			0,
			& \text{otherwise}.
		\end{cases}
	\end{aligned}
	\]
	Thus all mixed free cumulants of
	$\widetilde X_1,\ldots,\widetilde X_n$
	vanish, which proves that these
	operators are freely independent.
	Moreover, 
	$\widetilde X_j\sim\MP(\lambda_j,\alpha),
	 1\leq j\leq n.$	
	Therefore, the families
	\begin{equation}\label{eq:joint}
	    (X_1,\ldots,X_n) \text{ in $(\A,\varphi)$}
	\quad\text{and}\quad
	(\widetilde X_1,\ldots,\widetilde X_n)
\text{ in $(\CC,\tau_p)$}
	\end{equation}
	have the same joint $\ast$-distribution. 
	Since $\sum_{j=1}^n E_j=1_{\DD}$,
	\[
	\begin{aligned}
		\sum_{j=1}^n\widetilde X_j
		&=
		Z^{\frac1{2}}
		\left(\sum_{j=1}^n E_j\right)
		Z^{\frac1{2}}=Z.
	\end{aligned}
	\]
	Thus, we may identify $F$ with $Z$
	at the level of joint $*$-distributions. 
	Since $Z$ is invertible in the corner
	algebra $\CC$, we obtain
	\[
	\begin{aligned}
		Z^{-\frac1{2}}\widetilde X_jZ^{-\frac1{2}}
		&=
		Z^{-\frac1{2}}Z^{\frac1{2}}E_j
		Z^{\frac1{2}}Z^{-\frac1{2}}=pE_jp.
	\end{aligned}
	\]
Recall that  $p\in\operatorname{vN}(Z)$, and
	$\operatorname{vN}(Z)$ is free from
	$\DD$. 
	By the free-compression theorem
	\cite[Application~1.13]{NS96},
	$Z$ is free from 
	$\{pE_1p,\ldots,pE_np\}$   in $(\CC,\tau_p)$. 
	Finally, the equality \eqref{eq:joint} of the joint
	$\ast$-distributions, together with continuous
	functional calculus, implies that
	\[
	(F,B_1,\ldots,B_n) \text{ in $(\A,\varphi)$} \quad \text{and} \quad (Z,pE_1p,\ldots,pE_np) \text{ in $(\CC,\tau_p)$}
    \]
	have the same joint $*$-distribution. 
		Therefore, $F$ is free from $\{B_1,\ldots,B_n\}$ as claimed.
\end{proof}

\section{Alternative derivation of Marchenko-Pastur law}\label{sec:alternative}
Finally, we provide an alternative derivation of the marginal distributions of $\MM_1$ and $\MM_2$ in Theorem \ref{thm:main1} by combining a simple trick with known results from the literature. This argument does not establish their  freeness. As before, $(\A,\varphi)$ is a tracial $W^\ast$-probability space throughout this section. 
 We need the following facts. 

\begin{lemma}[{\cite[Corollary 7.2]{Capitaine}}]
\label{casalis} Let $X$ and $Y$ be freely independent and self-adjoint elements in $(\A,\varphi)$ with common distribution $\MP(1,1)$. Then $(X+Y)^{-1/2}X(X+Y)^{-1/2}$ has the arcsine distribution on $(0,1)$:
\begin{align}\label{freearcsine}
\frac{1}{\pi\sqrt{x(1-x)}}\one_{(0,1)}(x)\,dx.
\end{align}
\end{lemma}

\begin{lemma}[{\cite[Theorem 4.2]{Szpojankowski2015}}]
\label{szpojan}
 Let $X,Y$ be freely independent and self-adjoint elements in $(\A,\varphi)$ with distributions $\MP(\lambda,\alpha)$ and $\MP(\kappa,\alpha)$, respectively, where $\lambda,\kappa>0$, $\alpha>0$, and $\lambda+\kappa>1$. Then the elements  $(X+Y)^{-1/2}X(X+Y)^{-1/2}$
 and $X + Y$ are free.
\end{lemma}
\begin{lemma}\label{lemma:arcsine}
Let $X$ be a self-adjoint element in $(\A,\varphi)$ with the arcsine distribution
on $(0,1)$. 
Then $(2X-1_\A)^2$ has the same distribution as $X$.
\end{lemma}

\begin{proof}
Set $Y=2X-1_\A$. Then $Y$ has the symmetric arcsine
distribution on $(-1,1)$. By the characterization
of Arnold and Groeneveld \cite{ArnoldGroeneveld1980}, we get 
\[
Y^2\overset{d}{=}\frac{1_\A+Y}{2}=X. \qedhere
\]
\end{proof}

We use the $S$-transform; the reader is referred to e.g.~\cite{nica} for its definition and properties. If $X, Y$ are self-adjoint elements having the arcsine distribution \eqref{freearcsine} and Marchenko-Pastur distribution $\MP(\lambda,\alpha)$ with $\lambda,\alpha >0$, respectively, then their $S$-transforms are given by 
\begin{align}
S_{X}(z)=\frac{2+z}{1+z}, \qquad
 S_Y(z) = \frac{1}{\alpha(z+\lambda)};  \label{eq:StransormataBainomial}
\end{align}
the former formula can be deduced from the Cauchy transform 
\[
\int_0^1\frac{1}{z-x}\frac{dx}{\pi\sqrt{x(1-x)}} = \frac{1}{\sqrt{z(z-1)}}, 
\]
where $\sqrt{z(z-1)}$ denotes the analytic branch on
$\mathbb{C}\setminus[0,1]$ that is positive for real $z>1$.
The latter formula follows from the $R$-transform.

Now we show that
$ \MM_1\text{ and } \MM_2$ in Theorem \ref{thm:main1} 
 have the standard Marchenko-Pastur distribution. Suppose that $\SS$ and $\TT$ are freely independent standard semicircular elements. 
Using the operators 
\[
\XX:=\left(\frac{\SS+\TT}{\sqrt{2}}\right)^2\quad \text{and}\quad\YY:=\left(\frac{\SS-\TT}{\sqrt{2}}\right)^2,
\] 
 we can
 express $\MM_1$ in the form 
\begin{align}
\MM_1&=(\XX+\YY)^{-\frac1{2}}(\XX-\YY)^2(\XX+\YY)^{-\frac1{2}}.\label{eq:pom1}
 \intertext{Let $\FF:=\XX+\YY$ and $\GG:=(\XX+\YY)^{-\frac1{2}}\XX(\XX+\YY)^{-\frac1{2}}$. Then we can write   }
\MM_1&=\FF^{-\frac1{2}}(2\FF^{\frac1{2}}\GG\FF^{\frac1{2}}-\FF)^2\FF^{-\frac1{2}} \nonumber
\\&=(\FF-2\FF\GG-2\GG\FF+4\GG\FF\GG)\nonumber 
\\&=(2\GG-1_\A)\FF(2\GG-1_\A). \nonumber
\end{align}
By the orthogonal invariance of a standard semicircular family \cite[Theorem 5.3]{Nica96}, $(\SS+\TT)/\sqrt2$ and $(\SS-\TT)/\sqrt2$ are freely independent standard semicircular elements. 
Therefore $\XX$ and $\YY$ are free and both have distribution $\MP(1,1)$.
Lemma~\ref{szpojan} implies that $\FF$ and $\GG$ are free, and
$\FF\sim\MP(2,1)$. 
Since $(2G-1_\A)F(2G-1_\A)$ and $F(2G-1_\A)^2$ have the same moments, we have 
\begin{align*}
 S_{\MM_1}(z)&=S_{(2\GG-1_\A)\FF(2\GG-1_\A)}(z)=S_{(2\GG-1_\A)^2}(z)S_{\FF}(z). 
 \intertext{By Lemma~\ref{lemma:arcsine} and \eqref{eq:StransormataBainomial}, this equals}
 &S_{\GG}(z)S_{\FF}(z)=\frac{2+z}{1+z} \times \frac{1}{2+z}=\frac{1}{1+z}.
 \end{align*}
This is the $S$-transform of $\MP(1,1)$, so $\MM_1\sim\MP(1,1)$.
The same argument applies to $\MM_2$, which has the form~\eqref{eq:pom1}
with $\XX=\SS^2$ and $\YY=\TT^2$. These elements are free and both have
distribution $\MP(1,1)$; hence $\MM_2\sim\MP(1,1)$ as well.   

\section*{Acknowledgments}
The authors would like to thank Marek Bo\.zejko for suggesting research topics, many fruitful discussions,
and helpful comments. 
This research was funded in part by Narodowe Centrum
Nauki, Poland, under the WEAVE-UNISONO grant
2022/04/Y/ST1/00008 (Wiktor Ejsmont).

\section*{Statement on AI use}
The authors used ChatGPT to assist with translation and language editing of the manuscript. The authors reviewed and verified all AI-assisted content and take full responsibility for the scientific content of the paper.

\end{document}